\documentclass[12pt]{article}
\usepackage{mathrsfs}
\usepackage{latexsym,lineno}
\usepackage{epsfig}
\usepackage{color}
\usepackage{amsmath}\usepackage{fleqn}\usepackage{verbatim}\usepackage{epsf}
\usepackage{amsthm}\usepackage{graphicx, float}\usepackage{graphicx}
\usepackage{amsfonts}\usepackage{amssymb}\usepackage{graphpap}
\usepackage{epic}\usepackage{curves}
\newcommand{\be}{\begin{equation}}
\newcommand{\ee}{\end{equation}}
\newcommand{\benum}{\begin{enumerate}}
\newcommand{\eenum}{\end{enumerate}}
\newcommand{\bit}{\begin{itemize}}
\newcommand{\eit}{\end{itemize}}

\begin{document}
\def\s{\subseteq}
\def\n{\noindent}
\def\se{\setminus}
\def\dia{\diamondsuit}
\def\la{\langle}
\def\ra{\rangle}

%--------------------------------------------------------------

\title{ Total domination polynomials of graphs}

\footnotetext{The first and third authors are
partially supported by the Summer Graduate Research Assistantship Program of Graduate School, the second author  is partially supported  by the National Nature Science Foundation of China (Grant No.11171207). }
\author{ Jiuhua Hu$^a$, Erfang Shan$^b$, Shaohui Wang$^{a,c}$\footnote{  Corresponding authors: J. Hu (e-mail: jhu2@go.olemiss.edu), E. Shan (e-mail: efshan@shu.edu.cn), S. Wang (e-mail: swang4@go.olemiss.edu), C. Wang (wcxiang@mailccnu.edu.cn),  B. Wei (e-mail:  bwei@olemiss.edu). } , Chunxiang Wang$^d$, Bing Wei$^a$ \\
\small\emph {a. Department of Mathematics, The University of Mississippi, University, MS 38677, USA}\\
\small\emph {b. School of Management, Shanghai  University, Shanghai 200444, P.R.China}\\
\small\emph {c. Computer Science and Mathematics Department, Adelphi University
Garden City, NY 11530, USA}\\
\small\emph {d. School of Mathematics and Statistics, Central China Normal University
  Wuhan, 430079, P.R. China}}
\date{}
\maketitle

\begin{abstract}
Given a graph $G$,  a total dominating set $D_t$  is a vertex set that every vertex of $G$  is adjacent to some vertices of $D_t$ and let $d_t(G,i)$ be  the number of all total dominating sets with size $i$. The total domination polynomial, defined as $D_t(G,x)=\sum\limits_{i=1}^{| V(G)|} d_t(G,i)x^i$,  recently has been one of the considerable extended research in the field of domination theory. In this paper, we obtain the vertex-reduction and edge-reduction formulas of total domination polynomials. As  consequences, we give the total domination polynomials for paths and cycles. Additionally, we determine the sharp upper bounds of total domination polynomials for trees and  characterize the corresponding graphs attaining such bounds. Finally, we use the reduction-formulas to investigate  the
relations between vertex sets and total domination polynomials in $G$.

\end{abstract}

\section{Introduction}
Throughout this paper $G = (V, E)$ is a  finite simple undirected graph with vertex set $V = V(G)$ and edge set $E =E(G)$.  Let  $|V|$ denote the order of $G$. Especially, $G = \phi$ or $V = \phi$ if $|V| = 0$.
 For any $v\in V(G)$,  $N_G(v)=\{ w\in V(G):  vw\in E(G)\}$ is
the  open neighborhood of $v$ and $N_G[v]=N_G(v) \cup \{v\}$ is
the  closed neighborhood of $v$ in $G$.  If $N_G(v) = \phi$, then $v$ is called an isolated vertex.  The set of vertices adjacent in G to a vertex of a
vertex subset $S\subseteq V(G)$  is the open neighborhood   $N_G(S)$ of $S$,
 $N_G[S]=N_G(S)\cup S$ is the  closed neighborhood of $S$ and  $G-S$ is a subgraph induced by $V(G)-S$.
For $u \in V(G)$,
$G/ u$ is the contracted graph by the removal of $u$ and the addition of edges between any pair of non-adjacent neighbors of $u$,
$G \ominus u$ or $G\ominus u \ominus v$ represents a subgraph induced by $V(G) - N_G[u]$ or $V(G) - N_G[u] - N_G[v]$ respectively.  In particular, we set $ G-u\ominus v = (G-u)\ominus v $  and $ G-e\ominus u = (G-e)\ominus u$. 
   A vertex is said to be a pendant vertex if its open neighborhood contains exactly one vertex. The neighbor of a pendant vertex is called  a supporting vertex. A graph $F$ is called a forest if it has no cycles. When $F$ contains only one component, we say $F$ is a tree.
 For   graphs $G_1,G_2$,  $G_1\cup G_2$ is called the union graph of $G_1$ and $G_2$ with vertex set $V(G_1) \cup V(G_2)$ and edge set $E(G_1) \cup E(G_2)$, $G_1\vee G_2$ is called the  join  graph of $G_1$ and $G_2$ with vertex set $V(G_1)\cup V(G_2)$ and edge set $E(G_1)\cup E(G_2)\cup \{ uv: u\in V(G_1) \mbox{ and } v\in V(G_2)\}$. In [4,8],
 The $2$-corona of a graph $G$  is  defined by the graph of order $3| V(G)|$ obtained by attaching a path of length $2$ to each vertex of $G$ such that the resulting paths are vertex disjoint. Let $S_n$ , $T_n$ be the star, tree of order $n$.

A vertex set $D_t$ of a graph $G$ is a total dominating set \cite{cdh} if every vertex of $V(G)$  is adjacent to some vertices of $D_t$. Let $\gamma_t(G)$ be the minimum size of total dominating sets,  $\mathcal{D}$$_t(G,i)$ be the set of all total dominating sets of size $i$ and set $ d_t(G,i) = |$$\mathcal{D}$$_t(G,i)$$|$.
 The total domination polynomial \cite{2012}, defined as $D_t(G,x)=\sum\limits_{i=1}^{| V(G)|} d_t(G,i)x^i$, is one of   the extended research-area of the domination theory.   The domination polynomial was studied recently by several authors, see \cite{total1978, ap, survey1, survey, note, cdh, chromatic}.
 
The following equalities, which are easy to check by the concepts,  are very useful in calculating    total domination polynomials of graphs.
\vskip 2mm \noindent {\bf Proposition 1 } \emph{ Let $G$, $P_n$ and $C_n$ be a graph, a path and a cycle with $n$ vertices. Then\\
$(i)$ $D_t(G,x) \neq 0$ if and only if $G$ has no isolated vertices.\\
$(ii)$ The number of supporting vertices of $G$ is $n - d_t(G, n-1)$.}\\
$(iii)$\cite{bvv} $D_t(G,x) = D_t(G_1,x)D_t(G_2,x)$, $where$ $G= G_1 \cup G_2$.\\
$(iv)\;  D_t(P_1,x) = 0, D_t(P_2,x) = x^2, D_t(P_3,x) = x^3+2x^2,D_t(P_4,x) = x^4+2x^3+x^2. $\\
$(v)\cite{2012}\;  D_t(C_3,x) = x^3+3x^2, D_t(C_4,x) = x^4+4x^3+4x^2, D_t(C_5, x) = x^5+5x^4+5x^3, D_t(C_6,x) = x^6+6x^5+9x^4.$

The following proposition is due to
Cockayne et al.(1980)  and
Brigham et al.(2000).
\vskip 2mm \noindent {\bf Proposition 2 \cite{ hg, cdh}}
 \emph{ If $G$ is a connected graph of order $n \geq 3$, then
$2 \leq \gamma_t(G)\leq 2n/ 3$, the right equality holds if and only if $G$ is $C_3, C_6$ or  $2$-corona.
}

In general, it is very hard to find the total domination number(polynomial) of a graph and determine the solutions of related extremal problems. 
A plenty of properties of domination number are explored, see \cite{001,003,004,002,006}.
However, only a few classes of graphs with exact determination of the coefficients have been appeared in the literture. 
Vijayan and Kumar \cite{2012}(2012) obtained some properties  of total domination polynomials for cycles.
 Chaluvaraju et al.\cite{bvv}(2014)   presented some basic  properties of total domination polynomials  and  graph operations of the union and join of graphs.
However, the vertex and edge reduction formulas, which are very important tools to investigate the properties of graph polynomials, are still unknown for total domination polynomials.

This article not only obtains  recurrence relations of total domination polynomials for graph operations, but also digs out some interesting results on extremal problems and the characterization of  graphic structures  using the total domination polynomials of  a graph.
The main results of this study are detailed below.
\begin{enumerate}
 \item We present the vertex-reduction and edge-reduction formulas of total domination polynomials in Theorems 1  and 3. As   consequences, Theorems 2 and 4 give the recurrence relations for total domination polynomials of paths and cycles.
 \item We obtain the sharp upper bounds of total domination polynomials for trees and  characterize the corresponding graphs attaining such bounds in Theorem 5, which is a classic type of extremal problems.
 \item  Using the reduction-formulas, we investigate the values of total domination polynomials at $x = -1$ in Theorems 6 and 7. Our results indicate that  the number of all total dominating sets of even size and  that of odd size  can differ by at most 1 for the forest.
 \item  A direct relation between the vertices of degree 2 and the coefficients of total domination polynomials is given in Theorem 8.

 \end{enumerate}

\section{Reduction-formulas of total domination polynomials}
In this section, we determine the vertex-reduction and edge-reduction formulas of total domination polynomials.
For any graph $G$, let $W\subseteq V(G)$ and $C_W$ be a statement on $W$.  Denote $D_t(G,x)\{C_W\}$ to be the generating function for the number of total dominating sets of $G$ under  the condition $C_W$:
$$D_t(G,x)\{C_W\}=\sum_{ N_G[W]=V(G)}\varphi(W)x^{| W|}, \mbox{where } \varphi(W)= \left\{ \begin{array}{rcl}
1, && \mbox{if } C_W \mbox{ holds for }W,
  \\0,&&   \mbox{if otherwise}.
\end{array}\right.$$
For instance, if $w \in V(G)$, then $D_t(G,x)\{w \notin W\}$ represents the  total domination polynomial  generated by all total domination sets $W$ of $G$ with  $w \notin W$.
Also, define a new identity function $1_{G\ominus u \ominus v}$ as follows: For the vertices $u,v \in V(G)$,
$$1_{G\ominus u \ominus v} = \left\{ \begin{array}{rcl}
1,\;\;\;\;\;\;\;\;\;\;\;\;\;\;\;\;\;\;\;\;\;\;\;\; & \mbox{if } G\ominus u \ominus v = \phi,
  \\D_t(G \ominus u \ominus v,x), & \mbox{if }  G\ominus u \ominus v \neq \phi.
\end{array}\right.$$

\vskip 2mm \noindent {\bf Theorem 1} \emph{
For any connected graph $G$ and  $u\in V(G)$,
$$\begin{array}{rcl}
D_t(G,x) &=& D_t(G-u,x)+xD_t(G/u,x)-(1+x)D_t(G/ u,x)\{ N_G(u)\cap W =\emptyset\} \\&& +\sum_{v \in N_G(u)}x^21_{G\ominus u \ominus v}.
\end{array}$$
}
\noindent{\bf Proof.}  We first consider the case when the total dominating set $W$ does not contain $u$. Then
 $W$ is a total dominating set of $G$  if and only if $N_G(u) \cap W \neq \phi$ and $W$ is a total dominating set of $G-u$. Thus,
\begin{eqnarray}
D_t(G,x)\{u \notin W\} = D_t(G-u,x)\{N_G(u) \cap W \neq \phi\}.
\end{eqnarray}
Also,
when $W$ is a total domination set of $G-u$,  $W$ is also a total domination set of $G/u$ since $G-u$ is a subgraph by the definition of $G-u$ and $G/u$; Conversely,
if $W$ is a total domination set of $G/u$ and $N_G(u) \cap W = \phi$, then $W$ is still a total domination set of $G-u$. Thus,
\begin{eqnarray}
D_t(G-u,x)\{ N_G(u)\cap W=\phi\} = D_t(G/u,x)\{ N_G(u)\cap W=\phi\}.
\end{eqnarray}
 By (1) and (2), we can get that,
\begin{eqnarray}
D_t(G,x)\{u \notin W\}  &\stackrel{(1)}{=}&D_t(G-u,x)\{N_G(u) \cap W \neq \phi\} \nonumber \\
 &\stackrel{(2)}{=}& D_t(G-u,x)-D_t(G/u,x)\{N_G(u)\cap W=\phi\}.
\end{eqnarray}

Now we consider the case when   the total dominating set $W$ contains $u$. Then  $N_G(u)\cap W\neq \phi$.
If $|N_G(u)\cap W| = 1$, say $v \in N_G(u)\cap W$, and $| N_{G-u}(v)\cap W|=0$, then either  $G \ominus u\ominus v = \phi$ or  $W-\{u,v\}$ is a total dominating set of $G \ominus u\ominus v$.  We have
\begin{eqnarray}
&&D_t(G,x)\{u\in W,N_G(u)\cap W=\{ v\},| N_{G-u}(v)\cap W|=0\} \nonumber \\&&
= \left\{ \begin{array}{rcl}
x^2,\;\;\;\;\;\;\;\;\;\;\;\;\;\;\;\;\;\;\;\;\;\;\;\;\;\; & \mbox{if } G\ominus u \ominus v = \phi,
  \\x^2D_t(G \ominus u \ominus v,x), & \mbox{if }  G\ominus u \ominus v \neq \phi.
\end{array}\right.  =  x^21_{G\ominus u \ominus v}.
\end{eqnarray}
If $N_G(u)\cap W=\{ v\}$ and $| N_{G-u}(v)\cap W|\geq 1$, then $W-u$ is also a total dominating set of  $G/u$. Also, for any total dominating set $W'$ of $G/u$, we have $W'\cup\{u\}$ is a total dominating set of $G$. Thus,
\begin{eqnarray}
&&D_t(G,x)\{u\in W,N_G(u)\cap W=\{ v\},| N_{G-u}(v)\cap W|\geq 1\} \nonumber \\&&
=x D_t(G/u,x)\{N_G(u)\cap W=\{v\},| N_{G/u}(v)\cap W| \geq 1\}.
\end{eqnarray}
Similarly, if $| N_G(u)\cap W|\geq 2$, then $W-u$ is a total dominating set of $G/u$  and for any total dominating set $W'$ of $G/u$, then $W'\cup\{u\}$ is a total dominating set of $G$. Thus,
\begin{eqnarray}
&& D_t(G,x)\{u\in W,| N_G(u)\cap W|\geq 2\}=x D_t(G/u,x)\{| N_G(u)\cap W|\geq 2\}.
\end{eqnarray}
Furthermore, for the total dominating set $W$ of $G$,
if $v \in N_G(u)$ and $N_{G/u}(v)\cap W = \phi$, then $d_G(u) = d_G(v) =1$ and $uv \in E(G)$, that is, $v$ is an isolated vertex of $G/u$. Thus,
$W$ is not a total dominating set of $G/u$, that is,
 \begin{eqnarray} D_t(G/u,x)\{N_G(u) \cap W =\{v\}, |N_{G/u}(v) \cap W| = 0\}
 = 0.\end{eqnarray}
For $G/u$,     we can obtain that
 \begin{eqnarray}
&&D_t(G/u,x) - D_t(G-u,x)\{| N_G(u)\cap W|=0\} \nonumber \\&
\stackrel{(2)}{=}&D_t(G/u)\{|N_G(u)\cap W| \geq 1\}\nonumber
\\ &\stackrel{(7)}{=}& \sum_{v \in N_G(u)}D_t(G/u,x)\{N_G(u) \cap W  =   \{v\},    |N_{G/u}(v) \cap W| \geq 1\} \nonumber \\&&+D_t(G/u,x)\{|N_G(u) \cap W| \geq 2\}.
\end{eqnarray}
Also, by (4)-(8) and $v \in N_G(u)$,  one can obtain that
$$\begin{array}{rcl}
&& D_t(G,x)\{u\in W\} \\ &=&
D_t(G,x)\{u\in W,| N_G(u)\cap W|=1\}
+D_t(G,x)\{u\in W,| N_G(u)\cap W|\geq 2\} \\
&=&  \sum_{v \in N_G(u)} D_t(G,x)\{u\in W,N_G(u)\cap W=\{ v\}\}
+D_t(G,x)\{u\in W,| N_G(u)\cap W|\geq 2\} \\
&\stackrel{(7)}{=}& \sum_{v \in N_G(u)}D_t(G,x)\{u\in W,N_G(u)\cap W=\{ v\},| N_{G-u}(v)\cap W|\geq 1\} \\
&&+D_t(G,x)\{u\in W,| N_G(u)\cap W|\geq 2\}\\
&\stackrel{(4,5,6)}{=}&    \sum_{v \in N_G(u)} [x^21_{G\ominus u \ominus v} +  x D_t(G/u,x)\{N_G(u)\cap W=\{v\},| N_{G/u}(v)\cap W| \geq 1\}] \\ && + xD_t(G/u,x)\{| N_G(u)\cap W|\geq 2\} \\
&\stackrel{(8)}{=}& \sum_{v \in N_G(u)} x^21_{G\ominus u \ominus v}
+x [D_t(G/ u,x) -  D_t(G-u,x)\{| N_G(u)\cap W|=0\}].
\end{array}$$
Finally, combining (3) and the above equality, we can obtain that  $$\begin{array}{rcl}
D_t(G,x)&=&D_t(G,x)\{u\in W\}+D_t(G,x)\{u\notin W\} \\&=& D_t(G-u,x)+xD_t(G/u,x)-(1+x)D_t(G/ u,x)\{ N_G(u)\cap W =\emptyset\} \\&&+ \sum_{v \in N_G(u)} x^21_{G\ominus u \ominus v}
\end{array}$$  and Theorem 1 is true. $\hfill\Box$

By Theorem 1, we can get the following corollary.
\vskip 2mm \noindent {\bf Corollary 1.1} \emph{
For any connected graph $G$, if either $(i)$ there exist two vertices $u,v\in V(G)$ such that $N_G [v] \subseteq  N_G [u]$ or $(ii)$ there exists a supporting vertex $w \in N_G(u)$, then
$$
D_t(G,x)=D_t(G-u,x)+x D_t(G/ u,x)+\sum_{v\in N_G(u)}x^21_{G\ominus u \ominus v}.$$
}
\noindent {\bf Proof.}
Let $W\subseteq V(G)$  such that $N_G(u)\cap W=\phi$. For $(i)$, since $v$ cannot be dominated by $W$, then $D_t(G/u,x)\{N_G(u)\cap W=\phi\}=0$. For $(ii)$, we know that any supporting vertex $w$ must belong to any total dominating set $W'$ of $G$. Otherwise, the pendant vertices adjacent to $w$ can not be dominated by $W'$. Thus, $D_t(G/u,x)\{N_G(u)\cap W=\phi\}=0$.
Hence, by Theorem 1, we have Corollary 1.1 is true.$\hfill\Box$

Apply Corollary 1.1 to a path $P_n = v_1v_2...v_n$ with $n \geq 5$. Let $u = v_2, v \in N_{P_n}(u) = \{v_1,v_3\} $,  by Proposition $1(i)$ and $(iii)$,  we can get the recurrence relation of paths below.
\vskip 2mm \noindent {\bf Theorem 2} \emph{
Let $P_n$ be a path of order $n\geq 5$,
$$
D_t(P_{n},x)=xD_t(P_{n-1},x)+x^2D_t(P_{n-3},x)+x^2D_t(P_{n-4},x).
$$
}
As an immediate consequence, we give the total domination polynomials of paths as follows.
\vskip 2mm \noindent {\bf Corollary 2.1} \emph{ Let $P_n$ be a path of order $n \geq 1$,
$$D_t (P_n,x) = \left\{ \begin{array}{rcl}
\frac{2x^{\frac{n}{2}}}{x+4} + p(n),\;\;\;\;\;\;\;\;\;\;\;\;\;\;\;\;\;&\mbox{if }& n=4m,
  \\ -\frac{(x^2+3x)x^{\frac{n-1}{2}}}{x+4} + p(n),\;\;\;& \mbox{if }& n = 4m+1,\\
-\frac{2x^{\frac{n}{2}}}{x+4} + p(n), \;\;\;\;\;\;\;\;\;\;\;\;\;\;& \mbox{if}& n = 4m+2,\\
\frac{(x^2+3x)x^{\frac{n-1}{2}}}{x+4} + p(n),\;\;\;\; \;\;& \mbox{if}& n = 4m+3,
\end{array}\right.$$
}
where $p(n) = \frac{(x+2- \sqrt{x^2+4x})(x-\sqrt{x^2+4x})^n + (x+2 + \sqrt{x^2+4x})(x+\sqrt{x^2+4x})^n}{2^{n+1}(x+4)}. $

\noindent {\bf Proof.} If  $n =  1,2, 3, 4$, by   Proposition $1(iv)$,  Corollary 2.1 is true. If $n \geq 5$,
 the characteristic polynomial of the recursion of Theorem 2 is $\lambda^4-x\lambda^3-x^2\lambda-x^2 = 0$ with roots  $\lambda_1 =\sqrt{-x},
\lambda_2 =-\sqrt{-x},
\lambda_3 =\frac{x+\sqrt{x(x+4)}}{2},
\lambda_4 =\frac{x-\sqrt{x(x+4)}}{2}$. Hence,  $$D_t(P_n,x) = \alpha_1(x)\lambda_1^n+\alpha_2(x)\lambda_2^n+\alpha_3(x)\lambda_3^n+ \alpha_4(x)\lambda_4^n$$ with $\alpha_1(x),\alpha_2(x,)\alpha_3(x)$ and $\alpha_4(x)$ not dependent on $n$.
Using the initial conditions in Proposition $1(iv)$, we can form a system of equations  for $i =1,2,3,4 $,  $D_t(P_i,x) = \alpha_1(x)\lambda_1^i+\alpha_2(x)\lambda_2^i+\alpha_3(x)\lambda_3^i+ \alpha_4(x)\lambda_4^i$. By solving this system, we get that
$\alpha_1(x)= \frac{2+(x+3)\sqrt{-x}}{2(x+4)},
\alpha_2(x)= \frac{2-(x+3)\sqrt{-x}}{2(x+4)},
\alpha_3(x)= \frac{x+2+\sqrt{x(x+4)} }{2(x+4)},
\alpha_4(x)=\frac{x+2-\sqrt{x(x+4)} }{2(x+4)}.
$

Thus, by  $D_t(P_n,x) =\alpha_1(x)\lambda_1^n + \alpha_2(x)\lambda_2^n+ \alpha_3(x)\lambda_3^n+ \alpha_4(x)\lambda_4^n$ and setting $n = 4m+i$  with $i  = 0,1, 2,3$ respectively, we have
 Corollary 2.1 is true.
$\hfill\Box$

\vskip 2mm \noindent {\bf Theorem 3} \emph{
For any connected graph $G$ and  $e=uv\in E(G)$,
$$
D_t(G,x)=D_t(G-e,x)+x^2 1_{G\ominus u \ominus v}
+(1+x)[D_t(G-e\ominus u,x)\{v\in W\}+D_t(G-e \ominus v,x)\{u\in W\}].$$
}
\noindent{\bf Proof.}
Let $W$ be a total dominating set of $G$, we will consider the cases that whether $u,v \in W$ or not as follows.

\noindent{\bf Case 1. $\{u,v\}\subseteq W$.}

If $| N_G(u)\cap W|\geq 2 \mbox{ and } | N_G(v)\cap W|\geq 2$, then $W$ is also a total dominating set of $G-e$;
If $| N_G(u)\cap W|=1 \mbox{ and } | N_G(v)\cap W|=1$, that is, $N_G(u)\cap W=\{v\}$ and $N_G(v)\cap W=\{u\}$, then $W-\{u,v\}$ is a total dominating set of $G \ominus u \ominus v$ unless $G \ominus u \ominus v=\phi$;
If $| N_G(u)\cap W|=1 \mbox{ and } | N_G(v)\cap W|\geq 2$, then $W-\{u\}$ is a total dominating set of $G-e\ominus u$;
If $| N_G(u)\cap W|\geq 2 \mbox{ and } | N_G(v)\cap W|=1$, then $W-\{v\}$ is a total dominating set of $G-e\ominus v$.

Hence,
if $ G \ominus u \ominus v \neq \phi$, then
$D_t(G,x)\{\{u,v\}\subseteq W\}=D_t(G-e,x)\{\{u,v\}\subseteq W\}+x^2 D_t(G \ominus u \ominus v, x)
+x D_t(G-e\ominus u,x)\{v\in W\}+x D_t(G-e\ominus v,x)\{u\in W\};$
If $ G \ominus u \ominus v= \phi$, then
$D_t(G,x)\{\{u,v\}\subseteq W\}=D_t(G-e,x)\{\{u,v\}\subseteq W\}+x^2
+x D_t(G-e\ominus u,x)\{v\in W\}+x D_t(G-e\ominus v,x)\{u\in W\}.$

\noindent{\bf Case 2. $\{u,v\}\cap W=\phi$.}

Then  $W$ is also a total dominating set of $G-e$, that is,
$D_t(G,x)\{\{u,v\}\cap W=\phi\}=D_t(G-e,x)\{\{u,v\}\cap W=\phi\}.$

\noindent {\bf Case 3. $u\in W, v\notin W. $}

Since $W$ is a total dominating set of $G$,   $N_G(v)\cap W\neq\phi$.
If $| N_G(v)\cap W|\geq 2$, then $ W$ is also a total dominating set of $G-e$;
If $N_G(v)\cap W=\{u\}$, then $W$ is a total dominating set of $G-e\ominus v$.
Hence, we have
$D_t(G,x)\{u\in W,v\notin W\}=D_t(G-e,x)\{u\in W,v\notin W, W\cap N_{G-e}(v)\neq \phi\}
+ D_t(G-e\ominus v,x)\{u\in W\}.$

\noindent{\bf Case 4. $u\notin W, v\in W $.}

Similar to Case 3, we have
$D_t(G,x)\{u\notin W,v\in W\}=D_t(G-e,x)\{v\in W,u\notin W,W\cap N_{G-e}(u)\neq \phi\}+D_t(G-e\ominus u,x)\{v\in W\}.$

Finally, we have $$\begin{array}{rcl}
 D_t(G,x) &=&
D_t(G,x)\{\{u,v\}\subseteq W\}+D_t(G,x)\{\{u,v\}\cap W=\phi\} \\&& +D_t(G,x)\{u\in W,v\notin W\}  +D_t(G,x)\{u\notin W,v\in W\} \\
&=&  D_t(G-e,x)+x^2 1_{G\ominus u \ominus v}
+(1+x)[D_t(G-e\ominus u,x)\{v\in W\}\\&&+D_t(G-e \ominus v,x)\{u\in W\}]
\end{array}$$ and Theorem 3 is true.
$\hfill\Box$

By applying  Theorems 2 and 3, we can obtain the recurrence relations of total domination polynomials for cycles.

\vskip 2mm \noindent{\bf Theorem 4} \emph{
For any cycle $C_n=v_1v_2\cdots v_nv_1$ with $n \geq 7$,
$$
D_t(C_n,x)=xD_t(C_{n-1},x)+x^2D_t(C_{n-3},x)+x^2D_t(C_{n-4},x).
$$
}
\noindent{\bf Proof.}
Let $P_n = C_n-e = v_1v_2...v_n,$ where $e=v_1v_n$. We begin with the proof by Claim 1.

{\small\bf Claim 1.} {\em
$D_t(P_{n},x )\{v_{n}\in W\}=xD_t(P_{n-1},x )\{v_{n-1}\in W\}+x^2D_t(P_{n-3},x )\{v_{n-3}\in W\}+x^2D_t(P_{n-4},x )\{v_{n-4}\in W\}.$
}

{\it Proof of Claim 1.} For $i \geq 1$,
let $S_0 = \{ W: W\in$ $\mathcal{D}$$_t(P_n,i)$$, v_{n}\in W \}$, $S_1 = \{ W_1\cup \{v_n\}: W_1\in $ $\mathcal{D}$$_t(P_{n-1},i-1)$$, v_{n-1}\in W_1\}$, $S_2 = \{ W_2\cup \{v_n,v_{n-1}\}: W_2\in $ $\mathcal{D}$$_t(P_{n-3},i-2)$$, v_{n-3}\in W_2 \}$ and $S_3 = \{ W_3\cup \{v_n,v_{n-1}\}: W_3\in $ $\mathcal{D}$$_t(P_{n-4},i-2)$$, v_{n-4}\in W_3 \}$. Now it is enough to show that $S_0 = S_1 \cup S_2 \cup S_3$.

It is obvious that $S_1 \cup S_2 \cup S_3 \subset S_0$. Conversely, choose  any $W \in S_0$. The definition of total dominating set yields that
 $v_{n-1}\in W$. If $v_{n-2}\in W$, set $W_1=W-\{v_n\}$, that is, $W_1 \in $ $\mathcal{D}$$_t(P_{n-1},i-1)$ and $ W_1 \cup \{v_n\} \in S_1$.
 If $v_{n-2}\notin W$ and $v_{n-3}\in W$, then  $v_{n-4}\in W$ and therefore set $W_2= W-\{v_n,v_{n-1}\}$, that is, $W_2 \in$ $\mathcal{D}$$_t(P_{n-3},i-2)$ and  $W_2 \cup \{v_n,v_{n-1}\} \in S_2$.
 If $v_{n-2}\notin W$ and $v_{n-3}\notin W$, then
 $v_{n-4}\in W, v_{n-5}\in W$ and therefore set $W_3=W - \{v_n,v_{n-1}\}$, that is, $W_3 \in $ $\mathcal{D}$$_t(P_{n-4},i-2)$ and  $W_3 \cup \{v_n,v_{n-1}\} \in S_3$. Thus,
$S_0 \subset S_1 \cup S_2 \cup S_3$.
Hence, $D_t(P_{n},x )\{v_{n}\in W\}=xD_t(P_{n-1},x )\{v_{n-1}\in W\}+x^2D_t(P_{n-3},x )\{v_{n-3}\in W\}
+x^2D_t(P_{n-4},x )\{v_{n-4}\in W\}$ and Claim 1 is true. $\hfill\Box$

By Theorems 2, 3  and Claim 1, we have
$$\begin{array}{rcl}
D_t(C_n,x)&\stackrel{\mbox{Theorem 3}}{=}&D_t(P_{n},x)+x^2D_t(P_{n-4},x)+(1+x)(D_t(C_n-e-v_1-v_2,x)\{v_n\in W\}\\
&&+ D_t(C_n-e-v_n-v_{n-1},x)\{v_1\in W\})\\
&=&D_t(P_n,x)+x^2D_t(P_{n-4},x)+2(1+x)D_t(P_{n-2},x)\{v_1\in W\}\\
&\stackrel{\mbox{Theorem 2}}{=}&
[xD_t(P_{n-1},x)+x^2D_t(P_{n-3},x)+x^2D_t(P_{n-4},x)]\\
&&+x^2[ xD_t(P_{n-5},x) +x^2D_t(P_{n-7},x)+x^2D_t(P_{n-8},x)]\\
&&+ 2(1+x)[xD_t(P_{n-3},x)\{v_1\in W\}+x^2D_t(P_{n-5},x)\{v_1\in W\}\\
&&+x^2D_t(P_{n-6},x)\{v_1\in W\}]\\
&\stackrel{\mbox{Simplify}}{=}&x(D_t(P_{n-1},x)+x^2D_t(P_{ n-5},x)+2(1+x)D_t(P_{n-3},x)\{v_1\in W\})\\
&&+x^2(D_t(P_{n-3},x)+x^2D_t(P_{ n-7},x)+2(1+x)D_t(P_{n-5},x)\{v_1\in W\})\\
&&+x^2(D_t(P_{n-4},x)+x^2D_t(P_{ n-8},x)+2(1+x)D_t(P_{n-6},x)\{v_1\in W\})\\
&\stackrel{\mbox{Claim 1}}{=}& xD_t(C_{n-1},x)+x^2D_t(C_{n-3},x)+x^2D_t(C_{n-4},x).
\end{array}$$
Thus, Theorem 4 is true.
$\hfill\Box$

As a direct consequence, we obtain the total domination polynomials of cycles below.
\vskip 2mm \noindent {\bf Corollary 4.1} \emph{ For any cycle $C_n$ with $ n\geq 3$ vertices,
$$D_t (C_n,x) = \left\{ \begin{array}{rcl}
2(-x)^{\frac{n}{2}} + q(n),&\mbox{if }& n=2m,
  \\ q(n),\;\;\;\;\;\;\;\;\;\;\;\;\;\;\; & \mbox{if }& n = 2m+1,
\end{array}\right.$$
where $q(n) = \frac{(x-\sqrt{x^2+4x})^n+(x+\sqrt{x^2+4x})^n}{2^{n}}$.
}

\noindent{\bf Proof.}
If $n = 3, 4,5$ or $6$, by Proposition $1(v)$,   Corollary 4.1 is true. If $n \geq 7$,  the characteristic polynomial of the recursion of Theorem 4 is $\lambda^4-x\lambda^3-x^2\lambda-x^2 = 0$ with roots  $\lambda_1 =\sqrt{-x},
\lambda_2 =-\sqrt{-x},
\lambda_3 =\frac{x+\sqrt{x(x+4)}}{2},
\lambda_4 =\frac{x-\sqrt{x(x+4)}}{2}$. Thus,
$$D_t(C_n,x) = \alpha_1(x)\lambda_1^n+\alpha_2(x)\lambda_2^n+\alpha_3(x)\lambda_3^n+ \alpha_4(x)\lambda_4^n$$
with $\alpha_1(x),\alpha_2(x,)\alpha_3(x)$ and $\alpha_4(x)$ not dependent on $n$. Using the initial conditions in Proposition $1(v)$, we form a system of  equations  for $i =3,4,5,6 $, $D_t(C_i,x) = \alpha_1(x)\lambda_1^i+\alpha_2(x)\lambda_2^i+\alpha_3(x)\lambda_3^i+ \alpha_4(x)\lambda_4^i$. By solving this system, we get
$\alpha_1(x)= \alpha_2(x)=\alpha_3(x)=\alpha_4(x)=1.$

Thus, by $D_t(C_n,x) =\lambda_1^n + \lambda_2^n+ \lambda_3^n+ \lambda_4^n$ and setting $n = 2m+i$ with $i = 0, 1 $, we have Corollary 4.1 is true.
$\hfill\Box$

\section{On total domination polynomials of special graphs}
In this section, we obtain the sharp upper bounds of the coefficients of total domination polynomials for trees and  characterize the corresponding graphs attaining such bounds.

\vskip 2mm \noindent {\bf Theorem 5} \emph{
Let $T_n$ be a tree  in  $\mathcal{T}$$_n$ that is a class of trees with $n$ vertices.\\
$(i)$ $d_t(T_n,i) \leq  {n-1 \choose i-1}$ and the equality holds if and only if $T_n= S_n$. \\
$(ii)$ There is no tree $T'$   in  $\mathcal{T}$$_n$ such that  $ d_t(T_n,i) \geq  d_t(T',i)$ for each $i \geq 2$.
}

\noindent {\bf Proof.}
If $n = 0, 1$, then $D_t(T_n, x) = 0$, that is , Theorem $5(i)$ is true. Next we only consider the case that $n \geq 2$.
By Proposition 2, $\gamma_t(T_n) \geq 2$ and
it is easy to see that
$d_t(S_n, i) =  {n-1 \choose i-1}$.
Let $r$ be the number of supporting  vertices of
$T_n$. Since $n \geq 2$, then $ r \geq 1$. By ${n-r \choose i-1} \leq {n-1 \choose i-1}$, we obtain that  $S_n$ attains the maximal value for all the coefficients.   Thus, $(i)$ is true.

Now we will prove $(ii)$ by the contradiction. Assume that $T_1$ attains the minimal value at all of the coefficients. By Proposition $1(ii)$, we have
$d_t(T_1, n-1) = n-r$. Since $d_t(T_1, n-1)$ attains the minimal value, then
$r$ is as big as possible. Thus,  $ r  = \lfloor{\frac{n}{2}}\rfloor$, that is, $d_t(T_1, n-1) = n - \lfloor{\frac{n}{2}}\rfloor = \lceil{\frac{n}{2}}\rceil $.
On the other hand,
since $D_t(G,x) = \sum_{i=\gamma_t(G)}^n d_t(G, i)x^i$ and in order to obtain that every coefficient is as small as possible, then we would need that $\gamma_t$ is as big as possible. By Proposition 2, choose a tree $T_2$ such that $\gamma_t(T_2) = \lfloor{\frac{2n}{3}}\rfloor$. Thus, $d_t(T_2, \lfloor{\frac{2n}{3}}\rfloor - 1) = 0$.

Since $T_1$ contains only two types of vertices: $V_{p}$ and $V_{s}$, where $V_p$ and $V_s$ are the sets of all pendant vertices and supporting vertices of $T_1$ respectively, then
every vertex of $V_p \cup V_s$ is dominated by some supporting vertices. Otherwise, $T_1$ is not a connected graph. Thus,  $\gamma_t(T_1) = \lfloor{\frac{n}{2}}\rfloor$ and $d_t(T_1, \lfloor{\frac{2n}{3}}\rfloor- 1)  >  0 = d_t(T_2, \lfloor{\frac{2n}{3}}\rfloor - 1)$, a contradiction with the choice of $T_1$. Therefore, Theorem $5(ii)$ is true. $\hfill\Box$

By setting $x = 1$, $D_t(G, 1)$ is the number of all total domination sets in $G$. Togethering with Theorem 5, we have the following corollary.
\vskip 2mm \noindent {\bf Corollary 5.1} \emph{ Let $\mathcal{T}$$_n$ be a class of trees with $n$ vertices, then $S_n$ contains maximal number of total domination sets in $\mathcal{T}$$_n$.
}

By Corollary 2.1 and the proof of  Theorem 5, we can obtain the special values of total domination polynomials on certain graphs when $x = -1$.
\vskip 2mm \noindent {\bf Theorem 6} \emph{ Let $S_n, P_n$ be  a star and a path with  $n$ vertices. Then \\
$(i)$ $D_t(S_n, -1) = 1$;\\
$(ii)$ $D_t(P_{n_0},-1)=D_t(P_{n_2},-1)=D_t(P_{n_3},-1)=D_t(P_{n_5},-1)=1, D_t(P_{n_1},-1)=D_t(P_{n_4},-1)\\= 0$, where $n_i = i(\mbox{mod}\; 6)$  and $i = 0,1,...,5$.
}

\noindent {\bf Proof.}
$(i)$ Since $\gamma_t(S_n) = 2$ and
$D_t(S_n, -1) = \sum_{i=2}^{n}d_t(S_n, -1)(-1)^i = \sum_{i=2}^n{n-1 \choose i-1}(-1)^i $, then let $j = i -1$, we have
$$\begin{array}{rcl}
D_t(S_n, -1) &=& \sum_{j=1}^{n-1}{n-1 \choose j}(-1)^{j+1} \\
&=&  (-1)[\sum_{j=1}^{n-1}{n-1 \choose j}(-1)^j]\\
&=& (-1)[\sum_{j=0}^{n-1}{n-1 \choose j}(-1)^j - {n-1 \choose 0}(-1)^0]\\
 &=& (-1)[(1-1)^{n-1} -1] = 1.
\end{array}$$
For $(ii)$, by  Corollary 2.1 and setting $x = -1$, we have $D_t(P_{n}, -1) =\frac{2+cos(\frac{2n\pi}{3})-\sqrt{3}sin(\frac{2n\pi}{3})}{3}$.
 Let $n_i = 6m+i$ with $m\geq 0$ and $i = 0, 1, ..., 5$, one can obtain that
$$\begin{array}{rcl}
 && D_t(P_{n_0}, -1) = D_t(P_{n_3}, -1) =D_t(P_{n_5}, -1) = \frac{2+cos(4m\pi)-\sqrt{3}sin(4m\pi)}{3}=1,\\
&&  D_t(P_{n_1}, -1) = D_t(P_{n_4}, -1) = \frac{4sin^2(2m\pi)}{3}=0,\\
&& D_t(P_{n_2}, -1) = \frac{2+\sqrt{3}cos(\frac{1}{6}-4m)\pi+cos\frac{4(1+3m)\pi}{3}}{3}=1.
\end{array}$$
Thus,  Theorem 6 is true.
$\hfill\Box$

From Theorem 6, we see that $D_t(S_n,-1) = 1$ and $0 \leq D_t(P_n,-1) \leq 1$. In general, Theorem 7 shows that  the number of all total dominating sets of even size and  that of odd size for any forest can differ by at most 1.

\vskip 2mm \noindent {\bf Theorem 7} \emph{ Let $F_n$ be any forest with $n$ vertices, then
 $D_t(F_n,-1) \in \{0, 1\}$
}

\noindent {\bf Proof.}
By Proposition $1(iii)$, we will only need to consider $F_n$ as a tree. If $n=0,1$, then $D_t(F_n, -1) =0$ and Theorem 7 is true.  If $n = 2,3$, then $F_n$ is a path. By Proposition $1(iv)$, we have $D_t(F_n, -1) =1$ and the proof is done.

Now we will prove it by induction on $n \geq 4$. For $n=4$, either $F_n$ is a path $P_4$ or a star $S_4$. Also,  by Theorem 6, we have $D_t(P_4, -1) =0$ and $D_t(S_4, -1) =1$.
Next assume that $D_t(F_s,-1) \in \{0,1\}$ with $s < n$ and consider $F_n$. Choose $u$ to be any pendant vertex of $F_n$ and set $v \in N(u)$. By Theorem 1,  we have
$$\begin{array}{rcl}
D_t(F_n, -1) &=& D_t(F_n - u, -1) - D_t(F_n/u, -1) + \sum_{v\in N_G(u)}1_{G \ominus u \ominus v}\\
 &=& D_t(F_n - u, -1) - D_t(F_n/u, -1) + \sum_{v\in N_G(u),G \ominus u \ominus v=\phi}1\\
&&+ \sum_{v\in N_G(u),G \ominus u \ominus v \neq\phi} D_t(G \ominus u \ominus v, -1).
\end{array}$$ Since $F_n$ is a tree and $u$ is a pendant vertex, then
 $F_n-u = F_n/u$, that is, $D_t(F_n-u, -1) = D_t(F_n/ u, -1)$. For $G \ominus u \ominus v = \phi$, we have  $D_t(F_n, -1) =1$;
For $G \ominus u \ominus v \neq\phi$, we can obtain that $G \ominus u \ominus v$ must have at least one components $B_{i\geq 1}$. If one of the components is an isolated vertex, by Proposition $1(i)$ and $(iii)$, one can obtain that $D_t(F_n, -1) =0$. If every component $B_{i}$ is a tree with at least two vertices, by the induction hypothesis, we have  $D_t(B_i, -1)\in \{0,1\}$. Thus, by Proposition $1(iii)$, we have $D_t(G \ominus u \ominus v, -1) = \prod_{i\geq 1}D_t(B_i, -1) \in \{0,1\}$ for all the cases. Therefore, Theorem 7 is true.
$\hfill\Box$

The following result gives
the relation between  vertex set of degree 2 and  coefficients of total domination polynomials.
\vskip 2mm \noindent {\bf Theorem 8} \emph{
If $G$ is a connected graph of order $n$, then the number of vertices of degree 2  is at least ${n\choose 2}-{r\choose 2}-r(n-r)- d_t(G,n-2)$, where $r$ is the number of supporting vertices in $G$.
}

\noindent {\bf Proof.}
Let $L_0= \{v: d(v) = 2 \mbox{ and none}  \mbox{ of}\; N(v)  \mbox{ is  a supporting vertex}\}$, that is, $|L_0|$ is at most  the number of vertices of degree 2. It is enough to show that   $d_t(G,n-2)={n\choose 2}-{r\choose 2}-r(n-r)-|L_0| $.
Now we use "a pair" to give another representation of $L_0$.   Set
$L = \{\{ v_1,v_2\}: \mbox{there exists} \; v \in V(G) \mbox{ such that } N_G(v)=\{v_1,v_2\} \mbox{ and}  \mbox{ none of } v_1,v_2 \mbox{ is} \mbox{ supporting vertex}\big\}$, then $|L| = |L_0|$.
Next, for any vertices $ v_1,v_2 $ of $V(G)$,  we define
$$\begin{array}{rcl}
 && L_1=\{ V(G)-\{v_1,v_2\}: \mbox{there exists }\; v\in V(G) \mbox{ such }   \mbox{ that }   N_G(v)=\{v_1,v_2\}\}, \\
&&  L_2=\{ V(G)- \{v_1,v_2\}: \mbox{ at least one of } v_1,v_2 \mbox{ is a}  \mbox{ supporting vertex} \},\\
&& T=\{V(G) - \{ v_1,v_2\}: V(G)- \{ v_1,v_2\} \mbox{ is  not }    \mbox{a total dominating set of } G\}.
\end{array}$$
Then $|L|=|L_1|-|L_2|$ and $|L_2| = {r\choose 2}+r (n-r)$.  Now we first show Claim 2.

{\small\bf Claim 2.} {\em$T=L_1\cup L_2$.}

{\it Proof of Claim 2.}
 Clearly, $L_1\cup L_2\subseteq T$. Next we will show that $ T\subseteq L_1\cup L_2$. Take $W=V(G)- \{ v_1,v_2\} \in T$ and since W is not a total dominating set of $G$, by the definition of total dominating set, there exists $u\in V(G)$ such that
 $N_G(u)\subseteq \{v_1,v_2\}.$ If $u\in W$, then $d_G(u)\leq 2.$
 When $d_G(u)= 2$, then $N_G(u)=\{v_1,v_2\}$, that is, $W\in L_1$.
 When $d_G(u)= 1$, then $u$ is a pendant vertex. Without loss of generality, set $N_G(u)=\{v_1\}$, that is, $v_1$ is a supporting vertex and  $W\in S_2$.
 We know $|W| = n -2$, if $u\notin W$, then $u$ is either $v_1$ or $v_2$.
Without loss of generality, say $u=v_1$. Since $k=0$ and $N(u) \cap W = \phi$, then $N_G(u)={v_2}$, that is, $v_2$ is a supporting vertex. Thus, $W\in L_2$ and Claim 2 is true.  $\hfill\Box$

By the definition of total dominating set and Claim 2, we have
$$\begin{array}{rcl}
 d_t(G,n-2) &=&
 {n\choose 2}-|L_1\cup L_2|
 ={n\choose 2}-(| L_1|+| L_2|-| L_1\cap L_2| )\\
&=&  {n\choose 2}-({r\choose 2}+r (n-r))-(| L_1\mid- | L_1\cap L_2| )\\
&=& {n\choose 2}-{r\choose 2}-r(n-r)-| L_1-L_2| \\
&=& {n\choose 2}-{r\choose 2}-r(n-r)-|L|.
\end{array}$$
Therefore, Theorem 8 is true.
 $\hfill\Box$

\end{document}